\documentclass[12pt,reqno]{amsart}
\usepackage{amssymb}
\usepackage{cases}
\usepackage{bbm}
\usepackage{mathabx}
\usepackage{mathrsfs}
\usepackage{graphicx}
\usepackage{amsfonts}
\usepackage{enumerate}

\usepackage{amssymb, amsmath, amsfonts, latexsym}
\usepackage{enumerate}

\newtheorem{thm}{Theorem}[section]

\newtheorem{lemma}[thm]{Lemma}
\newtheorem{prop}[thm]{Proposition}

\newtheorem{remarks}[thm]{Remark}

\theoremstyle{definition}
\newtheorem{defn}{Definition}[section]
 \theoremstyle{remark}

\newcommand{\ee}{\mathbb{E}}

\newcommand{\rr}{\mathbb{R}}
\newcommand{\pp}{\mathbb{P}}

\newcommand{\e}{\varepsilon}

\def\AA{\mathcal A}
\def\BB{\mathcal B}

\def\DD{\mathcal D}
\def\FF{\mathcal F}
\def\EE{\mathcal E}

\def\SS{\mathcal S}

\def\G{\mathbf G}

\def\sgn{\rm sgn}

\def\AA{\mathcal A}
\def\BB{\mathcal B}

\def\DD{\mathcal D}

\def\GG{\mathcal G}
\def\FF{\mathcal F}
\def\EE{\mathcal E}
\def\HH{\mathcal H}

\def\SS{\mathcal S}

\def\<{\langle}
\def\>{\rangle}

\def\beq{\begin{equation}}
\def\nneq{\end{equation}}

\def\bdef{\begin{defn}}
\def\ndef{\end{defn}}

\def\bthm{\begin{thm}}
\def\nthm{\end{thm}}

\def\bprop{\begin{prop}}
\def\nprop{\end{prop}}

\def\brmk{\begin{remarks}}
\def\nrmk{\end{remarks}}

\def\bexa{\begin{exa}}
\def\nexa{\end{exa}}

\def\blem{\begin{lem}}
\def\nlem{\end{lem}}

\def\bcor{\begin{cor}}
\def\ncor{\end{cor}}

\def\<{\langle}
\def\>{\rangle}

\date{}

\def\bexe{\begin{exe}}
\def\nexe{\end{exe}}

\def\bprf{\begin{proof}}
\def\nprf{\end{proof}}

\def\bdes{\begin{description}}
\def\ndes{\end{description}}

\title[Moderate deviations for a fractional stochastic heat equation]{Moderate deviations for a fractional  stochastic heat equation with spatially correlated noise}

\author{Yumeng Li }
\address{Yumeng Li \\Department of Statistics and Finance, University of Science and Technology of China, Hefei, Anhui Province, 230026, China.}
\email{liyumeng@mail.ustc.edu.cn}

\author{Ran Wang}
\address{Ran Wang \\School of Mathematical sciences, University of Science and Technology of China, Hefei, Anhui Province, 230026, China.}
\email{wangran@ustc.edu.cn}

\author{Nian Yao} 

\address{Nian Yao\\College of Mathematics and Computational
Science, Shenzhen University, 518060,  Shenzhen, Guangdong Province, China.}
\email{yaonian@szu.edu.cn}

\author{Shuguang Zhang}
\address{Shuguang Zhang \\Department of Statistics and Finance, University of Science and Technology of China, Hefei, Anhui Province, 230026, China.}
\email{sgzhang@ustc.edu.cn}

\date{}

\begin{document}
\maketitle

 \noindent {\bf Abstract:}
 In this paper, we study the Moderate Deviation Principle for a perturbed  stochastic heat equation in the whole space $\rr^d, d\ge1$. This equation is driven by a Gaussian noise, white in time and correlated in space, and the differential operator is a fractional derivative operator.  The weak convergence method plays an important role.
 \vskip0.3cm

 \vskip0.3cm


\section{Introduction}
\noindent Since the work of Freidlin and Wentzell \cite{FW}, the theory of small perturbation large deviations for stochastic (partial) differential equation has been extensively developed (see \cite{DPZ, DZ}).
The large deviation principle (LDP) for stochastic  reaction-diffusion equations  driven by the space-time white noise was first  obtained by  Freidlin   \cite{Fre} and later by  Sowers \cite{Sowers},  Chenal and Millet \cite{CM}, Carrai and R\"ockner  \cite{CR} and other authors.  An LDP for a stochastic heat equation driven by a Gaussian noise, white in time and correlated in space was proved by M\'{a}rquez-Carreras and Sarr\`{a} \cite{MS}.  Recently, El Mellali and Mellouk proved an LDP for a fractional stochastic heat equation  driven by a spatially correlated noise in \cite{EM}.

 Like the  large deviations, the moderate deviation problems arise in the theory of statistical  inference quite naturally.  The   moderate deviation principle (MDP) can provide us with the rate  of convergence and a useful method for constructing asymptotic confidence intervals, see \cite{Erm}, \cite{GZ}, \cite{IK}, \cite{MiaoShen}, \cite{Wil}  and references therein.
Results on the MDP for processes with independent increments were obtained in De Acosta \cite{DeA}, Ledoux \cite{Led} and so on. The study of the MDP estimates for other processes has been carried out as well, e.g., Wu \cite{Wu} for Markov processes,  Guillin  and Liptser \cite{GL} for diffusion processes,  Wang and Zhang \cite{WZ} for  stochastic reaction-diffusion equations in $\rr$, Budhiraja {\it et al.} \cite{BDG} for stochastic differential equations with jumps, and references therein.

\vskip0.3cm
    In this paper, we study the MDP for the   fractional stochastic heat equation in spatial dimension $\rr^d$ driven by a spatially correlated noise.
  In \cite{LWZ}, we studied  the MDP for a perturbed  stochastic heat equations defined on $[0,T]\times [0,1]^d$,   driven by a spatially correlated noise. In that paper,  the method is the exponential approximation theorem (see \cite[Theorem 4.2.13]{DZ}),  which needs some exponential estimates.  However, due to the lack of good regularity properties of the Green function for the fractional heat equation, it is difficult to get those exponential estimates.  Instead
of proving exponential estimates, we will use the weak convergence approach (see \cite{BDM}) in this paper.

\vskip0.3cm
Now, let us give the  fractional stochastic heat equation
\begin{equation}\label{SPDE}
    \begin{cases}   \frac{\partial u^\e}{\partial t}(t,x)= \DD_{\underline{\delta}}^{\underline{\alpha}}u^{\e}(t,x)+b(u^\e(t,x))+ \sqrt\e\sigma(u^\e(t,x))\dot{F}(t,x),\\
        u^\e(0,x)=0,
                        \end{cases} \end{equation}
 where $\e>0, (t,x)\in [0,T]\times\rr^d,d\ge1, \underline{\alpha}=(\alpha_1,\cdots,\alpha_d),\underline{\delta}=(\delta_1,\cdots,\delta_d)$ and we will
assume that $\alpha_i\in]0,2]\setminus \{1\}$ and $|\delta_i|\le \min\{\alpha_i,2-\delta_i\}, i=1,\cdots, d$, $\dot F$ is the ``formal" derivative of the Gaussian perturbation and $\DD_{\underline{\delta}}^{\underline \alpha}$ denotes a non-local fractional differential operator on $\rr^d$ defined by
$$
\DD_{\underline{\delta}}^{\underline \alpha}:=\sum_{i=1}^d D_{\delta_i}^{\alpha_i}.
$$
Here $D_{\delta_i}^{\alpha_i}$ denotes the fractional differential derivative with respect to (w.r.t.) the \textit{i}-th coordinate defined via its Fourier transform $\FF$ by
$$
\FF(\DD_{\underline{\delta}}^{\underline \alpha}\phi)(\xi)=-|\xi|^{\alpha_i}\exp\left(-\imath\delta_i\frac{\pi}{2}\sgn{\xi}\right)\FF(\phi)(\xi),
$$
with $\imath^2+1=0$. The noise $F(t,x)$ is a martingale measure in the sense of Walsh \cite{Walsh}  and  Dalang \cite{Dal}, which will be defined with details in the sequel. The coefficients $b$ and $\sigma:\rr\rightarrow \rr$ are given functions. From now on, we shall refer to Eq. (\ref{SPDE}) as $Eq_{\underline{\delta},\e}^{\underline{\alpha}}(d, b,\sigma)$. See Section 2 for details.

 \vskip 0.4cm
As the parameter $\e$ tends to zero, the solutions $u^\e $ of \eqref{SPDE} will tend to the solution of the deterministic    equation defined by
 \begin{equation}\label{PDE}
    \begin{cases}   \frac{\partial u^0}{\partial t}(t,x)= \DD_{\underline{\delta}}^{\underline{\alpha}}u^{0}(t,x)+b(u^0(t,x)),\\
        u^0(0,x)=0.
                        \end{cases} \end{equation}

\vskip 0.2cm
  In this paper we shall investigate deviations of $u^\e $ from the deterministic solution $u^0$, as $\e$ decreases to $0$,
that is,  the asymptotic behavior of the trajectories,
$$
\frac{1}{\sqrt \e \lambda(\e)}\left(u^\e -u^0\right)(t,x),\quad (t,x)\in[0,T]\times\rr^d,
$$
where $\lambda(\e)$ is some deviation scale, which strongly influences the asymptotic behavior.

 The case $\lambda(\e)=1/\sqrt\e$ provides some large deviations estimates. Under suitable assumptions,  El Mellali and Mellouk proved that the law of the solution $u^{\e}$  satisfies a large deviation principle on the H\"older space in \cite{EM}.

  If $\lambda(\e)$ is identically equal to $1$, we are in the domain of the central limit theorem.

  To fill in the gap between the central limit theorem scale [$\lambda(\e)=1$] and the large deviations scale [$\lambda(\e)=1/\sqrt\e$],
we will study moderate deviations, that is when the deviation scale satisfies
\beq \label{h}
 \lambda(\e)\to+\infty,\ \ \sqrt\e \lambda(\e)\to0\ \ \ \text{as}\ \e\to0.
 \nneq
 The moderate deviation principle   enables us to refine the estimates obtained through the  central limit theorem. It provides the asymptotic  behavior  for  $\pp(\|u^\e-u^0\| \ge \delta \sqrt \e \lambda(\e))$ while the central limit theorem gives asymptotic bounds for
$\pp(\|u^\e-u^0\| \ge \delta \sqrt \e)$.
Throughout this paper, we assume \eqref{h} is in place.

\vskip0.3cm

The rest of this paper is organized as follows. In Section 2, the precise framework is stated.
In Section 3, the skeleton equation is studied. It is proved that the solution is a continuous map from the level set into the H\"older space.  Section 4  is devoted to the proof
of the moderate deviation principle by the weak convergence approach.
We give some precise estimates of the fundamental solution $G$ in the appendix.
\vskip 0.3cm

Throughout the paper, $C_p$ is a positive constant depending on the  parameter $p$, and $C, C_1,\cdots$ are   constants depending on no  specific parameter (except $T$ and the Lipschitz constants), whose value  may be different from line to line by convention.

 For any $T>0, K\subset\rr^d, \beta=(\beta_1,\beta_2)$, let $C^{\beta}([0,T]\times K;\rr^d)$  be the H\"older space equipped with the norm defined by
\beq\label{eq norm}
\|f\|_{\beta,K}:=\sup_{(t,x)\in[0,t]\times K}|f(t,x)|+\sup_{s\neq t\in [0,T]}\sup_{x\neq y\in K}\frac{|f(t,x)-f(s,y)|}{|t-s|^{\beta_1}+|x-y|^{\beta_2}}.
\nneq

Since $C^{\beta}([0,T]\times K; \rr^d)$ is not separable, we consider the space $C^{\beta',0}([0,T]\times K;\rr^d)$ of H\"older continuous $f$ with the degree $\beta'_i<\beta_i, i=1,2$ such that
$$
\lim_{\delta\rightarrow0^+}\left(\sup_{|t-s|+|x-y|<\delta}\frac{|f(t,x)-f(s,y)|}{|t-s|^{\beta_1'}+|x-y|^{\beta_2'}} \right)=0
$$
and $C^{\beta',0}([0,T]\times K;\rr^d)$ is a Polish space containing $C^{\beta}([0,T]\times K; \rr^d)$.
From now on, let $\EE^{\beta}([0,T]\times K; \rr^d):=C^{\beta,0}([0,T]\times K; \rr^d)$, where $\beta=(\beta_1,\beta_2)$.

\section{Framework}
In this section, let us give the framework taken from Boulanba {\it et al}. \cite{BEM} and  El Mellali and Mellouk \cite{EM}.

\subsection{The operator $\DD_{\underline{\delta}}^{\underline\alpha}$}
In one dimension space, the operator $D_{\delta}^{\alpha}$ is a closed, densely defined operator on $L^2(\rr)$ and it is the infinitesimal generator of a semigroup which is in general not symmetric and not a contraction. It is self-adjoint only when $\delta=0$ and in this case, it coincides with the fractional power of the Laplacian.

According to \cite{DI, koma}, $D_{\delta}^{\alpha}$ can be represented for $1<\alpha<2$, by
$$
D_{\delta}^{\alpha}=\int_{-\infty}^{+\infty}\frac{\phi(x+y)-\phi(x)-y\phi'(x)}{|y|^{1+\alpha}}(\kappa_-^{\delta}\mathbf{1}_{(-\infty,0)}(y)+\kappa_+^{\delta}\mathbf{1}_{(0,+\infty)}(y))dy,
$$
and for $0<\alpha<1$, by
$$
D_{\delta}^{\alpha}=\int_{-\infty}^{+\infty}\frac{\phi(x+y)-\phi(x)}{|y|^{1+\alpha}}(\kappa_-^{\delta}\mathbf{1}_{(-\infty,0)}(y)+\kappa_+^{\delta}\mathbf{1}_{(0,+\infty)}(y))dy,
$$
where $\kappa_-^{\delta}$ and $\kappa_+^{\delta}$ are two non-negative constants satisfying $\kappa_-^\delta+\kappa_+^\delta>0$ and $\phi$ is a smooth function for which the integral exists, and $\phi'$ stands for its derivative. This representation identifies it as the infinitesimal generator for a non-symmetric $\alpha$-stable L\'evy process.

Let $G_{\alpha,\delta}(t,x)$ denotes the fundamental solution of the equation $Eq_{\delta,1}^{\alpha}(1,0,0)$, that is, the unique solution of the Cauchy problem
 \begin{equation*}
    \begin{cases}   \frac{\partial u}{\partial t}(t,x)= D_{\delta}^{\alpha}u^{0}(t,x),\\
        u(0,x)=\delta_0(x), \ \ \ \ t>0, x\in\rr,
                        \end{cases} \end{equation*}
where $\delta_0$ is the Dirac distribution. Using Fourier's calculus one gets
\beq\label{eq fundamental solution}
G_{\alpha,\delta}(t,x)=\frac1{2\pi}\int_{-\infty}^{\infty}\exp\left(-\imath zx-t|z|^{\alpha}\exp\left(-\imath\delta\frac{\pi}{2}\sgn(z)\right)\right)dz.
\nneq
The relevant parameters, $\alpha$ called the index of {\it stablity} and $\delta$ called the {\it skewness }, are real numbers satisfying $\alpha\in ]0,2]$ and $|\delta|\le \min\{\alpha,2-\alpha\}$.

Now, for higher dimension $d\ge1$ and any multi index $\underline{\alpha}=(\alpha_1,\cdots\alpha_d)$ and $\underline{\delta}=(\delta_1,\cdots,\delta_d)$, let ${\mathbf G}_{\underline{\alpha},\underline{\delta}}(t,x)$ be the Green function of the deterministic equation $Eq_{\underline{\delta}}^{\underline{\alpha}}(d,0,0)$. Clearly,
\begin{align}
{\mathbf G}_{\underline{\alpha},\underline{\delta}}(t,x)=&\Pi_{i=1}^{d}G_{\alpha_i,\delta_i}(t,x_i)\notag\\
=&\frac{1}{(2\pi)^{d}}\int_{\rr^d}exp\left(-\imath \langle\xi,x\rangle-t\sum_{i=1}^d|\xi_i|^{\alpha_i}\exp\left(-\imath\delta_i\frac{\pi}{2}\sgn(\xi_i)\right)\right)d\xi,
\end{align}
where $\langle\cdot,\cdot\rangle$ stands for the inner product in $\rr^d$.

The properties of the Green function ${\mathbf G}_{\underline{\alpha},\underline{\delta}}(t,x)$ will be given in the appendix.

\subsection{The driving noise $F$ }
Let us explicitly describe here the spatially homogeneous noise, see Dalang \cite{Dal}.

Let $\SS(\rr^{d+1})$ be the space of Schwartz test functions. On a complete probability space $(\Omega, \GG,\pp)$, the noise $F=\{F(\phi),\phi\in\SS(\rr^{d+1})\}$ is assumed to be an $L^2(\Omega,\GG,\pp)$-valued Gaussian process with mean zero and covariance functional given by
$$
J(\varphi,\psi):=\ee\left[F(\phi)F(\psi)\right]=\int_{\rr_+}ds\int_{\rr^d}\left(\phi(s,\star)\ast\tilde \psi(s,\star)\right)(x)\Gamma(dx)ds,\ \ \ \phi,\psi\in\SS(\rr^{d+1}),
$$
where $\tilde\psi(s,x):=\psi(s,-x)$ and $\Gamma$ is a non-negative and non-negative definite tempered measure, therefore symmetric. The symbols $\ast$ denotes the convolution product and $\star$ stands for the spatial variable.

Let $\mu$ be the spectral measure of $\Gamma$, which is also a trivial tempered measure, that is $\mu=\FF^{-1}(\Gamma)$ and this gives
\beq\label{eq spectral measure}
J(\phi,\psi)=\int_{\rr_+}ds\int_{\rr^d}\mu(d\xi)\FF\phi(s,\star)(\xi)\overline{\FF\psi(s,\star)}(\xi),
\nneq
where $\bar z$ is the complex conjugate of $z$.

As in Dalang \cite{Dal}, the Gaussian process $F$ can be extended to a worthy martingale measure, in the sense of Walsh \cite{Walsh},
$$
M:=\{M_t(A),\ \ t\in \rr_+, A\in\BB_b(\rr^d)\},
$$
where $\BB_b(\rr^d)$ denotes the collection of all bounded Borel measurable sets in $\rr^d$.
Let $\GG_t$ be the completion of the $\sigma$-field generated by the random variables $\{F(s,A);0\le s\le t, A\in \BB_b(\rr^d)\}$.

Then,  Boulanba {\it et al.} \cite{BEM} gave a rigorous meaning to the solution of equation $Eq_{\underline{\delta},\e}^{\underline \alpha}(d,b, \sigma)$ by means of a joint measurable and $\GG_t$-adapted process $\{u^\e(t,x);(t,x)\in\rr_+\times\rr^d\}$ satisfying, for each   $t\ge0$ and for almost all $x\in\rr^d$ the following evolution equation
\begin{align}\label{u e}
 u^\e(t,x)=& \sqrt\e\int_0^t\int_{\rr^d}\G_{\underline{\alpha}, \underline{\delta}}(t-s,x-y)\sigma(u^\e(s,y))F(dsdy)\notag\\
 &+\int_0^tds\int_{\rr^d}\G_{\underline{\alpha}, \underline{\delta}}(t-s,x-y)b(u^\e(s,y))dy.
 \end{align}

In order to prove our main result, we are going to give other equivalent approach to the solution of $Eq_{\underline{\delta},\e}^{\underline \alpha}(d,b, \sigma)$, see \cite{DQ}.  To start with, let us denote by $\HH$ the Hilbert space obtained by the completion of $\SS(\rr^d)$ with the inner production
\begin{align*}
\langle \phi,\psi\rangle_{\HH}:=&\int_{\rr^d}\Gamma(dx)(\phi\ast \tilde \psi)(x)\\
=&\int_{\rr^d}\mu(d\xi)\FF\phi(\xi)\overline{\FF \psi}(\xi), \ \ \ \ \phi,\psi\in \SS(\rr^d).
\end{align*}
The norm induced by $\langle\cdot,\cdot\rangle_{\HH}$ is denoted by $\|\cdot\|_{\HH}$.

By the Walsh's theory of the martingale measures \cite{Walsh}, for $t\ge0$ and $h\in\HH$ the stochastic integral
$$
B_t(h)=\int_0^t\int_{\rr^d}h(y)F(ds,dy),
$$
is well-defined and the process $\{B_t(h); t\ge0, h\in\HH\}$ is a cylindrical Wiener process on $\HH$, that is:
\begin{itemize}                                                                                                              \item[(a)] for every $h\in\HH$ with $\|h\|_{\HH}=1$, $\{B_t(h)\}_{t\ge0}$ is a standard Wiener process,
 \item[(b)] for every $t\ge0, a,b\in \rr$ and $f,g\in\HH$, $$B_t(af+bg)=aB_t(f)+bB_t(g) \  \ \ \text{almost surely}.$$
 \end{itemize}

Let $\{e_k\}_{k\ge1}$ be a complete orthonormal system (CONS) of the Hilbert space $\HH$, then
$$
\left\{B_t^k:=\int_0^t\int_{\rr^d}e_k(y)F(ds,dy);k\ge1\right\}
$$
defines a sequence of independent standard Wiener processes and we have the following representation
\beq\label{eq wiener}
B_t:=\sum_{k\ge1}B_t^ke_k.
\nneq
Let $\{\FF_t\}_{t\in[0,T]}$ be the $\sigma$-field generated by the random variables $\{B_s^k;s\in[0,t], k\ge1\}$. We define the predictable $\sigma$-field in $\Omega\times[0,T]$ generated by the sets $\{(s,t]\times A;A\in\FF_s, 0\le s\le t\le T\}$. In the following, we can define the stochastic integral with respect to cylindrical Wiener process $(B_t(h))_{t\ge0}$ (see e.g. \cite{DPZ} or \cite{DQ}) of any predictable square-integrable process with values in $\HH$ as follows
$$
\int_0^t\int_{\rr^d}g\cdot dB:=\sum_{k\ge1}\int_0^t\langle g(s),e_k\rangle_{\HH} dB_s^k.
$$
  Note that the above series converges in $L^2(\Omega, \FF,\pp)$ and the sum does not depend on
the selected { CONS}. Moreover, each summand, in the above series, is a classical It\^o
integral with respect to a standard Brownian motion, and the resulting stochastic
integral is a real-valued random variable.

In the sequel, we shall consider the mild solution to equation  $Eq_{\underline{\delta},\e}^{\underline \alpha}(d,b, \sigma)$ given by
\begin{align}\label{SPDE solution}
u^{\e}(t,x)=&\sqrt\e\sum_{k\ge1}\int_0^t\langle \G_{\underline\alpha, \underline\delta}(t-s,x-\cdot)\sigma(u^\e(s,\star)), e_k\rangle_{\HH} dB_s^k\notag\\
&+\int_0^t\left[\G_{\underline\alpha,\underline \delta}(t-s)\ast b(u^\e(s,\star)) \right](x)ds,
\end{align}
for any $t\in[0,T],x\in\rr^d$.

\subsection{Existence, uniqueness and H\"older regularity to equation}
 For a multi-index $\underline{\alpha}=(\alpha_1,\cdots,\alpha_d)$ such that $\alpha_i\in]0,2]\setminus \{1\}, i=1,\cdots, d$ and any $\xi\in\rr^d$, let
 $$
 S_{\underline{\alpha}}(\xi)=\sum_{i=1}^d|\xi_i|^{\alpha_i}.
 $$
Assume the following assumptions on the functions $\sigma, b$ and the measure $\mu$:{\it
\begin{itemize}
  \item[{\bf (C)}:] The functions $\sigma$ and $b$ are Lipschitz, that is there exists some constant $L$ such that
\beq\label{Lip}
\|\sigma(x)-\sigma(y)\|\le L|x-y|, \ \ |b(x)-b(y)|\le L|x-y|\
\nneq
for all $ x,y\in\rr^d$.
  \item[(\bf $H_{\eta}^{\underline \alpha}$):] Let $\underline{\alpha}$ as defined above  and $\eta\in ]0,1]$, it holds that
 $$ \int_{\rr^d}\frac{\mu(d\xi)}{(1+  S_{\underline{\alpha}}(\xi))^{\eta}}<+\infty.$$
\end{itemize}}
The last assumption stands for an integrability condition w.r.t. the spectral measure $\mu$. Indeed, the following stochastic integral
$$
\int_0^T\int_{\rr^d}\G_{\underline{\alpha},\underline{\delta}}(T-s,x-y) F(ds,dy)
$$
is well-defined if and only if
$$
\int_0^Tds\int_{\rr^d}\mu(d\xi)|\FF\G_{\underline{\alpha},\underline{\delta}}(s,\ast)(\xi)|^2<+\infty.
$$
More precisely,  by \cite[Lemma 1.2]{BEM}, there exist two positive constants $c_1,c_2$ such that
\beq\label{eq H equivalent}
c_1\int_{\rr^d}\frac{\mu(d\xi)}{1+  S_{\underline{\alpha}}(\xi)}\le
\int_0^Tds\int_{\rr^d}\mu(d\xi)\left|\FF\G_{\underline{\alpha},\underline{\delta}}(s,\ast)(\xi)\right|^2 \le c_2 \int_{\rr^d}\frac{\mu(d\xi)}{1+  S_{\underline{\alpha}}(\xi)}.
\nneq
\vskip0.3cm

Under the assumptions $(C)$ and $(H_{\eta}^{\underline \alpha})$,   Boulanba {\it et al.} proved that Eq.  \eqref{SPDE solution}  admits a unique solution $u^\e$   such  that
\beq\label{eq solution u e}
\sup_{t\in[0,T]}\sup_{x\in\rr^d}\ee|u^\e(t,x)|^p<+\infty,\ \  \forall  T>0, p\ge2.
\nneq
See \cite[Theorem 2.1]{BEM}. Moreover, Theorem 3.1 in \cite{BEM} tells us that  the trajectories of the  solution $u^\e(t,x):(t,x)\in\rr_+\times\rr^d$ to Eq. \eqref{SPDE solution} are $\beta=(\beta_1,\beta_2)$-H\"older continuous in $(t,x)\in[0,T]\times K$ for every $K$ compact subset of $\rr^d$ and every $\beta_1\in (0,{1-\eta}/{2})$, $\beta_2\in (0,\min\{\alpha_0(1-\eta), 1/2\})$, where $\alpha_0:=\min_{1\le i\le d}\{\alpha_i\}$.

Consequently, the random field solution $\{u^\e(t,x);(t,x)\in[0,T]\times K\}$ to Eq. \eqref{SPDE solution} lives in the H\"older space $C^{\beta}([0,T]\times K;\rr^d)$ equipped with the norm $\|f\|_{\beta,K}$ given in \eqref{eq norm}.

Particularly, taking $\e=0$, the deterministic solution $u^0$ to \eqref{PDE}  has the following estimates
\begin{align}\label{eq u 0}
\sup_{t\in[0,T]}\sup_{x\in\rr^d}\ee|u^0(t,x)|^p<+\infty,\ \  \forall  T>0, p\ge2.
\end{align}
and $\|u^0\|_{\beta,K}<\infty$ for any compact set $K\subset\rr^d$.

\section{Skeleton equations}
The purpose of this section is to study the skeleton equation, which will be used in the weak convergence approach.
\vskip0.3cm
From now on, we furthermore suppose that{\it
  \begin{itemize}
        \item[({\bf D}):]\  The function $b$ is differentiable, and its derivative $b'$ is  Lipschitz. More precisely, there exists a positive constant $L'$ such that
\beq\label{H3}
|b'(y)-b'(z)|\le L'|y-z|\  \text{for all } y,z\in\rr.\ \
\nneq
      \end{itemize}}
Combined with the Lipschitz continuity of $b$, we  conclude that
\beq\label{H3'}
|b'(z)|\le L,\ \ \ \ \  \forall  z\in\rr.
\nneq

For $T>0$, let $\HH_T:=L^2([0,T]; \HH)$, which is a real separable Hilbert space such that, if $\varphi, \psi \in\HH_T$,
$$
\langle \varphi,\psi\rangle_{\HH_T}:=\int_0^T\langle \varphi(s,\cdot),\psi(s,\cdot)\rangle_{\HH}ds.
$$
Denote $\|\cdot\|_{\HH_T}$ the norm induced by $\langle \cdot,\cdot\rangle_{\HH_T}$. For any $N>0$, define
 $$
\HH_T^N:=\{h\in \HH_T;\|h\|_{\HH_T}\le N\},
$$
and we consider that $\HH_T^N$ is endowed with the weak topology of $\HH_T$.

For any $h\in \HH_T$, consider the deterministic evolution  equation (called Skeleton equation)
\begin{align}\label{eq Z h}
Z^h(t,x)=&\int_0^t\left\langle \G_{\underline\alpha, \underline \delta}(t-s, x-\star)\sigma(u^0(s,\star)),h(s,\star)\right\rangle_{\HH}ds\notag\\
&+\int_0^t\left[\G_{\underline\alpha, \underline \delta}(t-s)\ast ( b'(u^0(s,\star)Z^h(s,\star) )\right](x)ds,
\end{align}
where the first term on the right-hand side of the above equation can be written as
$$
\sum_{k\ge1}\int_0^t\left\langle  \G_{\underline\alpha, \underline \delta}(t-s, x-\star)\sigma(u^0(s,\star),e_k\right\rangle_{\HH}h_k(s)ds,
$$
with $h_k(t):=\langle h(t), e_k\rangle_{\HH}, t\in [0,T],k\ge1$.

\vskip0.3cm

 Using the strategy in the proof of Proposition 2.7 in \cite{EM}, one can obtain   the following result. Here we omit its  proof.
\begin{prop}\label{Prop skeleton} Assuming  conditions $(C), (H_{\eta}^{\underline{\alpha}})$ and $(D)$, there exists a unique solution $Z^h$ to Eq. \eqref{eq Z h}, which satisfies
\beq\label{eq Z h bound}
\sup_{h\in \HH_T^N}\sup_{(t,x)\in[0,T]\times\rr^d}|Z^h(t,x)|<+\infty.
\nneq
\end{prop}

\bthm\label{thm Skeleton}   Assuming  conditions $(C), (H_{\eta}^{\underline{\alpha}})$ and $(D)$, the mapping $h:\HH_T^N\rightarrow Z^h\in\EE^\beta([0,T]\times K; \rr^d)$ is continuous with respect to the weak topology, where $K$ is a compact set in $\rr^d$, $\beta=(\beta_1,\beta_2)$ satisfies that  $0<\beta_1<\alpha_0(1-\eta)/2, 0<\beta_2<1-\eta$ and $\alpha_0=\min_{1\le i\le d}\{\alpha_i\}$.

\nthm
\bprf Let $0<\beta_1<\alpha_0(1-\eta)/2, 0<\beta_2<1-\eta$ and  $\{h, (h_n)_{n\ge1}\}\subset \HH_T^N$ such that for any $g\in\HH_T$,
$$
\lim_{n\rightarrow\infty}\langle h_n-h, g \rangle_{\HH_T}=0.
$$
We need to prove that
\beq\label{eq convergence}
\lim_{n\rightarrow \infty}\|Z^{h_n}-Z^h\|_{\beta,K}=0.
\nneq
According to Lemma \ref{Lem 3} below (a particular case of Lemma A.1 in \cite{BMS}), the proof of \eqref{eq convergence} can be divided into two steps:
\begin{itemize}
  \item[(1)]{\bf Pointwise convergence}: for any $(t,x)\in [0,T]\times K$,
  \beq\label{eq pointwise}
  \lim_{n\rightarrow\infty}\left|Z^{h_n}(t,x)-Z^h(t,x)  \right| = 0.
  \nneq
  \item[(2)]{\bf Estimation of the increments}: for any $(t,x), (s,y)\in [0,T]\times K$,
  \begin{align}\label{eq increments}
  &\sup_{n\ge1}\left|(Z^{h_n}(t,x)-Z^h(t,x))-(Z^{h_n}(s,y)-Z^h(s,y))   \right|\notag\\
  \le& C\left(|t-s|^{\beta_1}+|x-y|^{\beta_2} \right).
  \end{align}
\end{itemize}
We will prove those two estimates in the following two steps.
\vskip0.3cm
{\bf Step 1}. {\it Pointwise convergence}.  For any $(t,x)\in [0,T]\times K$,
\begin{align}\label{eq Point 1}
&Z^{h_n}(t,x)-Z^h(t,x)\notag\\
=&\int_0^t\left\langle \G_{\underline\alpha, \underline \delta}(t-s, x-\star)\sigma(u^0(s,\star)),h_n(s,\star)-h(s,\star)\right\rangle_{\HH}ds\notag\\
&+\int_0^t\left\{\G_{\underline\alpha, \underline \delta}(t-s)\ast\left[  b'(u^0(s,\star)\left(Z^{h_n}(s,\star)-Z^h(s,\star)\right)\right]\right\}(x)ds\notag\\
=:&I_1^n(t,x)+I_2^n(t,x).
\end{align}

Since $h_n, h\in \HH_T^N$ and $u^0$ is bounded in $[0,T]\times\rr^d$, by Cauchy-Schwarz's inequality on the Hilbert space $\HH_T$, $({\bf C})$ and \eqref{eq H equivalent}, we have
\begin{align*}
|I_1^n(t,x)|^2\le &\int_0^t\|\G_{\underline\alpha, \underline \delta}(t-s, x-\star)\sigma(u^0(s,\star))\|_{\HH}^2ds\cdot\int_0^t\|h_n(s)-h(s)\|_{\HH}^2ds\notag\\
\le &4N^2C\int_0^t\|\G_{\underline\alpha, \underline \delta}(t-s, x-\star)\|_{\HH}^2ds\notag\\
\le & C(N)<+\infty,
\end{align*}
here $C(N)$ is independent of $n,t, x$. Since $h_n\rightarrow h$ weakly in $\HH_T^N$, we know that $I_1^n\rightarrow 0$ in $C([0,T]\times \rr^d;\rr)$ by Arz\`ela-Ascoli Theorem. This implies that
\beq\label{eq I 1}
\lim_{n\rightarrow \infty}\sup_{t\in[0,T],x\in\rr^d}|I_1^n(t,x)|=0.
\nneq

Set $\zeta^n(t):=\sup_{0\le s\le t, x\in\rr^d}|Z^{h_n}(s,x)-Z^h(s,x)|$. By Lemma \ref{lem fundamental solution} and \eqref{H3'}, we have
\begin{align}\label{eq I 2}
|I_2(t,x)|
\le &\int_0^t\int_{\rr^d}\G_{\underline\alpha, \underline \delta}(t-s, x-y)\left| b'(u^0(s,y)\left(Z^{h_n}(s,y)-Z^h(s,y)\right)\right|dyds\notag\\
\le & L\int_0^t\int_{\rr^d}\G_{\underline\alpha, \underline \delta}(t-s, x-y)\sup_{0\le l\le s,z\in\rr^d}\left| Z^{h_n}(s,z)-Z^h(s,z)\right|dyds\notag\\
\le& L\int_0^t \zeta^n(s)ds.
\end{align}

By \eqref{eq Point 1} and \eqref{eq I 2}, we have
$$
\zeta^n(t)\le L\int_0^t \zeta^n(s)ds+\sup_{t\in[0,T],x\in\rr^d}|I_1^n(t,x)|.
$$
Hence, by the Gronwall's lemma and \eqref{eq I 1}, we obtain that
$$ \zeta^n(T)\le e^{LT}\lim_{n\rightarrow\infty}\sup_{t\in[0,T],x\in\rr^d}|I_1^n(t,x)|\longrightarrow 0,\  \text{as } n\rightarrow\infty,
$$
which is stronger than \eqref{eq pointwise}.
\vskip0.3cm

{\bf Step 2}. {\it Estimation of the increments.} For any $0\le t \le T, s>0, x\in\rr^d,y\in K$,
\begin{align}\label{eq Z step 2}
&\left[Z^{h_n}(t+s,x+y)-Z^h(t+s,x+y)\right]-\left[Z^{h_n}(t,x)-Z^h(t,x)\right] \notag\\
= & \int_0^{t+s}\left\langle \G_{\underline\alpha, \underline \delta}(t+s-l, x+y-\star)\sigma(u^0(l,\star)),h_n(l,\star)-h(l,\star)\right\rangle_{\HH}dl\notag\\
&+\int_0^{t+s}\left\{\G_{\underline\alpha, \underline \delta}(t+s-l)\ast\left[  b'(u^0(l,\star)\left(Z^{h_n}(l,\star)-Z^h(l,\star)\right)\right]\right\}(x+y)dl\notag\\
&-\int_0^t\left\langle \G_{\underline\alpha, \underline \delta}(t-l, x-\star)\sigma(u^0(l,\star)),h_n(l,\star)-h(l,\star)\right\rangle_{\HH}dl\notag\\
&-\int_0^t\left\{\G_{\underline\alpha, \underline \delta}(t-l)\ast\left[  b'(u^0(l,\star)\left(Z^{h_n}(l,\star)-Z^h(l,\star)\right)\right]\right\}(x)dl\notag\\
= & \int_0^t\left\langle [\G_{\underline\alpha, \underline \delta}(t+s-l, x+y-\star)-\G_{\underline\alpha, \underline \delta}(t+s-l, x-\star)]\sigma(u^0(l,\star)),h_n(l,\star)-h(l,\star)\right\rangle_{\HH}dl\notag\\
& +\int_0^t\left\langle [\G_{\underline\alpha, \underline \delta}(t+s-l, x-\star)-\G_{\underline\alpha, \underline \delta}(t-l, x-\star)]\sigma(u^0(l,\star)),h_n(l,\star)-h(l,\star)\right\rangle_{\HH}dl\notag\\
&+\int_t^{t+s}\left\langle \G_{\underline\alpha, \underline \delta}(t+s-l, x+y-\star)\sigma(u^0(l,\star)),h_n(l,\star)-h(l,\star)\right\rangle_{\HH}dl\notag\\
&+\Bigg[\int_0^{t+s}\left\{\G_{\underline\alpha, \underline \delta}(t+s-l)\ast\left[  b'(u^0(l,\star)\left(Z^{h_n}(l,\star)-Z^h(l,\star)\right)\right]\right\}(x+y)dl\notag\\
&\ \ \ \  \ -\int_0^t\left\{\G_{\underline\alpha, \underline \delta}(t-l)\ast\left[  b'(u^0(l,\star)\left(Z^{h_n}(l,\star)-Z^h(l,\star)\right)\right]\right\}(x)dl\Bigg]\notag\\
=:&A_1^n+A_2^n+A_3^n+A_4^n.
\end{align}

Since $u^0$ is bounded, $h_n,h\in \HH_T^N$, we know that
\beq\label{eq finite}
\sup_{n\ge1}\int_0^T\|u^0(s,\star)(h^n-h)(s,\star)\|^2_{\HH}ds<\infty.
\nneq
 This, together with Cauchy-Schwarz's  inequality, \eqref{eq finite} and  Lemma \ref{lem holder},  implies that for each $0<\beta_1< (1-\eta)/2, 0<\beta_2<\min\{(1-\eta)\alpha_0/2, 1/2\}$, there exists a constant $C$ independent of $n$ such that
\beq\label{eq A1}
|A_1^n|\le C |y|^{\beta_2},\ \ |A_2^n|\le C s^{\beta_1},\ \
|A_3^n|\le C s^{\beta_1}.
\nneq

Let us now give the estimate of $A_4^n$.
Denote
$$
V_n(t,x):=b'(u^0(t,x))\left(Z^{h_n}(t,x)-Z^h(t,x)\right).
$$
By the estimates in Step 1, \eqref{H3'} and \eqref{eq Z h bound}, we know that
\beq\label{eq V n bound} \sup_{n\ge1}\sup_{t\in[0,T],x\in\rr^d}|V_n(t,x)|<\infty.\nneq

After a change of variable,  we have
\begin{align*}
A_4^n=&\int_0^s\int_{\rr^d}\G_{\underline{\alpha},\underline{\delta}}(t+s-l,x+y-z)V(l,z)dldz\\
&+\int_0^t\int_{\rr^d}\G_{\underline{\alpha},\underline{\delta}}(t-l,x-z)\left[V(s+l,y+z)-V(l,z) \right]dldz.
\end{align*}

By (i) of Lemma \ref{lem fundamental solution} and \eqref{eq V n bound}, we know that
$$
\sup_{n\ge1}\int_0^s\int_{\rr^d}\G_{\underline{\alpha},\underline{\delta}}(t+s-l,x+y-z)\left|V_n(l,z)\right|dldz\le c_1 s.
$$

By the Lipschitz continuity of $b'$,  the H\"older continuity of $u^0$ and the boundness of $Z^{h_n}$ and $Z^h$, we have
\begin{align}\label{eq V -}
&|V(s+l,y+z)-V(l,z)|\notag\\
\le & \left|b'(u^0(s+l, y+z))-b'(u^0(l,z))\right|\cdot\left|Z^{h_n}(s+l,y+z)-Z^{h}(s+l,y+z)\right|\notag\\
&+\left|b'(u^0(l,z))\right|\cdot\left|[Z^{h_n}(s+l,y+z)-Z^{h}(s+l,y+z)]-[Z^{h_n}(l,z)-Z^{h}(l,z)]\right|\notag\\
\le & C(L') \left|u^0(s+l,y+z)-u^0(l,z)\right|\notag\\
&+\left|b'(u^0(l,z))\right|\cdot\left|[Z^{h_n}(s+l,y+z)-Z^{h}(s+l,y+z)]-[Z^{h_n}(l,z)-Z^{h}(l,z)]\right|\notag\\
\le &c_1(s^{\beta_1}+|y|^{\beta_2})+c_2\sup_{z\in\rr^d} \left|[Z^{h_n}(s+l,y+z)-Z^{h}(s+l,y+z)]-[Z^{h_n}(l,z)-Z^{h}(l,z)]\right|.
\end{align}

Putting
$$\phi(l,s,y) :=\sup_{z\in\rr^d}\Big\{\left[Z^{h_n}(s+l,y+z)-Z^{h}(s+l,y+z)\right]-\left[Z^{h_n}(l,z)-Z^{h}(l,z)\right]\Big\}.$$
 Hence, by \eqref{eq Z step 2}, \eqref{eq A1} and  \eqref{eq V -}, we have
$$
\phi(t,s,y)\le c_3(s+s^{\beta_1}+|y|^{\beta_2})+c_2\int_0^t \phi(l,s,y)dl.
 $$
Therefore by the Gronwall's inequality,
\begin{align*}
&\sup_{0\le t\le T}\sup_{x\in\rr^d}\left|[Z^{h_n}(t+s,x+y)-Z^{h}(t+s,x+y)]-[Z^{h_n}(t,x)-Z^{h}(t,x)]\right|\\\notag
\le & c(s+s^{\beta_1}+|y|^{\beta_2}).
\end{align*}
The proof is complete.
\nprf

 \section{Moderate deviation principle }

 The main aim of this paper is to prove that $\frac{1}{\sqrt\e \lambda(\e)}(u^\e-u^0)$ satisfies an LDP on  the H\"older space, where $\lambda(\e)$ satisfies \eqref{h}.
This special type of LDP is usually called the moderate deviation principle of $u^\e$ (cf. \cite{DZ}).
\subsection{Large Deviation Principle}\
First,  recall the  definition of  large deviation principle. See \cite{DZ}.

Let $(\Omega,\mathcal{F},\mathbb{P})$ be a probability space with an increasing family $\{\FF_t\}_{0\le t\le T}$ of the sub-$\sigma$-fields of $\FF$ satisfying the usual conditions.
Let $\mathcal{E}$ be a Polish space with the Borel $\sigma$-field $\mathcal{B}(\mathcal{E})$.

    \bdef\label{def-Rate function}
         A function $I: \mathcal{E}\rightarrow[0,\infty]$ is called a rate function on
       $\mathcal{E}$,
       if for each $M<\infty$, the level set $\{x\in\mathcal{E}:I(x)\leq M\}$ is a compact subset of $\mathcal{E}$. A family of positive numbers $\{\lambda(\e)\}_{\e>0}$ is called a speed function if $\lambda(\e)\rightarrow +\infty$ as $\e\rightarrow 0$.
         \ndef

    \bdef\label{def LDP}  A family
       $\{X^\e\}$ of $\EE$-valued random elements is  said to satisfy the large deviation principle on $\mathcal{E}$
       with rate function $I$ and with speed function $\{\lambda(\e)\}_{\e>0}$, if the following two conditions
       hold.
       \begin{itemize}
         \item[$(a)$](Large deviation upper bound) For each closed subset $F$ of $\mathcal{E}$,
              $$
                \limsup_{\e\rightarrow 0}\frac1{\lambda(\e)}\log\mathbb{P}(X^\e\in F)\leq- \inf_{x\in F}I(x).
              $$
         \item[$(b)$](Large deviation lower bound) For each open subset $G$ of $\mathcal{E}$,
              $$
                \liminf_{\e\rightarrow 0}\frac1{\lambda(\e)}\log\mathbb{P}(X^\e\in G)\geq- \inf_{x\in G}I(x).
              $$
       \end{itemize}
    \ndef


\subsection{The main result}

From now on, we always assume that  $K$ is a compact set in $\rr^d$, $\underline\alpha=(\alpha_1,\cdots,\alpha_d),\alpha_i\in[0,2]\setminus\{1\}$ for $i=1,\cdots,d$, $\beta=(\beta_1,\beta_2)$ satisfies that  $0<\beta_1<\alpha_0(1-\eta)/2, 0<\beta_2<1-\eta$ where $\alpha_0=\min_{1\le i\le d}\{\alpha_i\}$, $\eta\in]0,1]$.
\vskip0.3cm
The main result of this paper is the following theorem.
\bthm\label{MDP}   Assuming  conditions $(C), (H_{\eta}^{\underline{\alpha}})$ and $(D)$ for $\eta\in]0,1]$. Let $u^\e$ be the solution of Eq. \eqref{SPDE solution}.  Then  the law of $(u^\e-u^0)/(\sqrt\e \lambda(\e))$ obeys an LDP on the space
$\EE^\beta([0,T]\times K; \rr^d)$ with speed $\lambda^2(\e)$ and with rate function
\beq\label{eq rate}
I(f)=\inf_{h\in \HH_T; Z^h=f}\left\{\frac12\|h\|_{\HH_T}^2\right\},
\nneq
where $Z^h$ is defined by \eqref{eq Z h}.
\nthm

\subsection{Proof of Theorem \ref{MDP}} \
\vskip0.3cm
\noindent
 We shall apply the weak convergence approach to establish moderate deviation principle.

\vskip0.3cm
Denote by $\AA_{\HH}$ the set of predictable process which belongs to $L^2(\Omega\times [0,T];\HH)$. For any $N>0$, let
$$
\AA_{\HH}^N:=\{h\in \HH_T; \|h\|_{\HH_T}\le N\}.
$$

For any $v\in\AA_{\HH}^N$ and $\e\in(0,1]$, define the controlled equation $Z^{\e,v}$ by
\begin{align}\label{eq Z e}
&Z^{\e,v}(t,x)\notag\\
=&\frac{1}{\lambda(\e)}\sum_{k\ge1}\int_0^t\left\langle \G_{\underline\alpha, \underline\delta}(t-s,x-\star)\sigma\big(u^0(s,\star)+\sqrt\e\lambda(\e)Z^{\e,v}(s,\star)\big), e_k  \right\rangle_{\HH}dB_s^k\notag\\
&+\int_0^t\left\langle \G_{\underline\alpha, \underline\delta}(t-s,x-\star)\sigma\big(u^0(s,\star)+\sqrt\e\lambda(\e)Z^{\e,v}(s,\star)\big), v(s,\star)  \right\rangle_{\HH}ds\notag\\
&+\int_0^t\left\{\G_{\underline\alpha, \underline\delta}(t-s)\ast\left[\frac{b\big(u^0(s,\star)+\sqrt\e\lambda(\e)Z^{\e,v}(s,\star)\big) -b(u^0(s,\star) )}{\sqrt\e\lambda(\e)} \right] \right\}ds.
\end{align}

Following the proof of Theorem 2.1 in \cite{BEM} and Proposition 2.7 in \cite{EM}, one can show the existence and uniqueness of the stochastic controlled equation given by \eqref{eq Z e}.
\begin{lemma}\label{lem control}
Assuming  conditions $(C), (H_{\eta}^{\underline{\alpha}})$ and $(D)$ for $\eta\in]0,1]$, there exists a unique random field solution to Eq. \eqref{eq Z e}, $\{Z^{\e,v}(t,x);(t,x)\in[0,T]\times\rr^d\}$, which satisfies that for any $p\ge1$,
\beq\label{eq Z e estimate}
\sup_{\e\le 1}\sup_{v\in\AA_{\HH}^N}\sup_{(t,x)\in[0,T]\times\rr^d}\ee[|Z^{\e,v}(t,x)|^p]<\infty.
\nneq
\end{lemma}

 Inspiriting by \cite{BDM}, let us consider the following two conditions which correspond on the weak convergence approach frame in our setting. Also refer to the weak convergence approach used in \cite{EM}.
 \begin{itemize}
   \item[(a)] The set $\{Z^h;h\in \HH_T^N\}$ is a compact set of $\EE^{\beta}$, where $Z^h$ is the solution of Eq. \eqref{eq Z h}.
   \item[(b)] For any family $\{v^\e;\e>0\}\subset\AA_{\HH}^N$ which converges in distribution as $\e\rightarrow 0$ to $v\in \AA_{\HH}^N$, as $\HH_T^N$-valued random variables, we have
       $$
       \lim_{\e\rightarrow 0}Z^{\e,v^{\e}}=Z^v\ \
       $$
       in distribution, as $\EE^{\beta}$-valued random variables, where $Z^v$ denotes the solution of Eq. \eqref{eq Z h} corresponding to the $\HH_T^N$-valued random variable $v$ (instead of a deterministic function $h$).
 \end{itemize}

\vskip0.3cm
\bprf[The proof of Theorem \ref{MDP}]
Applying to Theorem 6 in \cite{BDM},   a verification of  conditions (a) and (b) implies the validity of Theorem \ref{MDP}.
  Condition (a)  follows from the continuity of the mapping  $h:\HH_T^N\rightarrow Z^h\in\EE^\beta([0,T]\times K; \rr^d)$, which has been   established in Theorem \ref{thm Skeleton}. Next, we verify Condition (b).

By the Skorokhod representation theorem, there exist a probability $(\bar\Omega,\bar\FF,(\bar\FF_t),\bar\pp)$, and, on this basis, a sequence of independent Brownian motions $\bar B=(\bar B_k)_{k\ge1}$ and also a family of $\bar \FF_t$-predictable processes $\{\bar v^\e;\e>0\}, \bar v$ belonging to $L^2(\bar\Omega\times [0,T];\HH)$ taking values on $\HH_T^N$, $\bar\pp$-a.s., such that the joint law of $(v^\e,v, B)$ under $\pp$ coincides with that of  $(\bar v^\e,\bar v, \bar B)$ under $\bar\pp$ and
$$
\lim_{\e\rightarrow0}\langle\bar v^\e-\bar v, g \rangle_{\HH_T}=0, \ \ \forall g\in\HH_T,  \bar \pp-a.s..  \ \
$$

Let $\bar Z^{\e,\bar v^\e}$ be the solution to a similar equation as \eqref{eq Z e} replacing $v$  by $\bar v^\e$ and $B$ by $\bar B$. Thus, to verify the condition (b), it is sufficient to prove that
\beq\label{eq b conv}
\lim_{\e\rightarrow0} \|\bar Z^{\e,\bar v^\e}-\bar Z^{\bar v} \|_{\beta,K}=0 \text{  in probability} .
\nneq

From now on, we drop the bars in the notation for the sake of simplicity, and we denote
$$
Y^{\e,v^\e,v}:=Z^{\e,v^\e}-Z^v.
$$

According to Lemma \ref{Lem 3}, the proof of \eqref{eq b conv} can be divided into two steps: for any $ (t,x),(s,y)\in[0,T]\times K$ with $K$ being compact in $\rr^d$,
\begin{itemize}
  \item[(1)] {\bf Pointwise convergence}:
  \beq\label{eq Z 1}
  \lim_{\e\rightarrow0}|Y^{\e,v^\e,v}(t,x)|   =0 \ \ \ \text{in probability.}
    \nneq
  \item[(2)] {\bf Estimation of the increments}:  there exists a constant $C$ satisfied that
  \begin{align}\label{eq Z 2}
  \sup_{\e\le 1}\ee \left|Y^{\e,v^\e,v}(t,x)-Y^{\e,v^\e,v}(s,y)\right|^2\le C\left[|t-s|^{\beta_1}+|x-y|^{\beta_2}\right]
    \end{align}
\end{itemize}
Those two estimates will be established in  the following steps.

\vskip0.3cm
{\bf Step 1.} {\it Pointwise convergence.}  For any $(t,x)\in[0,T]\times K$,
\begin{align}\label{eq I0}
&Z^{\e,v^\e}(t,x)-Z^{v}(t,x)\notag\\
=&\frac{1}{\lambda(\e)}\sum_{k\ge1}\int_0^t\left\langle \G_{\underline\alpha, \underline\delta}(t-s,x-\star)\sigma\big(u^0(s,\star)+\sqrt\e\lambda(\e)Z^{\e,v^\e}(s,\star)\big), e_k  \right\rangle_{\HH}dB_s^k\notag\\
&+\Bigg\{\int_0^t\left\langle \G_{\underline\alpha, \underline\delta}(t-s,x-\star)\sigma\big(u^0(s,\star)+\sqrt\e\lambda(\e)Z^{\e,v^\e}(s,\star)\big), v^\e(s,\star)  \right\rangle_{\HH}ds\notag\\
&\ \ \ \ -\int_0^t\left\langle \G_{\underline\alpha, \underline\delta}(t-s,x-\star)\sigma\big(u^0(s,\star)\big), v(s,\star)  \right\rangle_{\HH}ds\notag\Bigg\}\\
&+\Bigg\{\int_0^t\left\{\G_{\underline\alpha, \underline\delta}(t-s)\ast\left[\frac{b\big(u^0(s,\star)+\sqrt\e\lambda(\e)Z^{\e,v^\e}(s,\star)\big) -b(u^0(s,\star) )}{\sqrt\e\lambda(\e)} \right](x) \right\}ds\notag\\
&\ \ \ \ -\int_0^t\left\{\G_{\underline\alpha, \underline\delta}(t-s)\ast\left[b'(u^0(s,\star))Z^v(s,\star)  \right](x) \right\}ds\Bigg\}\notag\\
=:&A_1^\e(t,x)+A_2^\e(t,x)+A_3^\e(t,x).
\end{align}

Let
\beq\label{eq J}
\mathcal J(t):=\int_{\rr^d}\mu(d\xi)|\FF\G_{\underline\alpha,\underline\delta}(t)(\xi)|^2.
\nneq

For the first term $A_1^\e$, by Burkholder's inequality,  the lipschitz property of $\sigma$, Eq. \eqref{eq u 0} and \eqref{eq Z e estimate}, we have that
 \begin{align}\label{eq A e 1}
\ee\left[|A_1^\e(t,x)|^2\right]=&\lambda^{-2}(\e)\ee\left[\int_0^t\left\|\G_{\underline\alpha,\underline\delta}(t-s,x-\star)\sigma(u^0(s,\star)+\sqrt\e\lambda(\e)Z^{\e,v}(s,\star))\right\|_{\HH}^2ds \right]\notag\\
\le &\lambda^{-2}(\e)c_1 \int_0^tds\left(1+\sup_{(r,y)\in[0,s]\times\rr^d}\ee|u^0(r,y)+\sqrt\e\lambda(\e)Z^{\e,v^\e}(s,y) |^2\right)\times \mathcal J(t-s)\notag\\
\le& c_2  \lambda^{-2}(\e).
\end{align}

The second term is further divided into two terms:
\begin{align*}
|A_2^\e(t,x)|\le &\int_0^t\left|\left\langle \G_{\underline\alpha, \underline\delta}(t-s,x-\star)\left[\sigma\big(u^0(s,\star)+\sqrt\e\lambda(\e)Z^{\e,v^\e}(s,\star)\big)-\sigma(u^0(s,\star))\right], v^\e(s,\star)  \right\rangle_{\HH}\right|ds\notag\\
&+\int_0^t\left|\left\langle \G_{\underline\alpha, \underline\delta}(t-s,x-\star)\sigma\big(u^0(s,\star)\big), v^\e(s,\star)-v(s,\star)  \right\rangle_{\HH}\right|ds\notag \\
=:&A_{2,1}^\e(t,x)+A_{2,2}^\e(t,x)
\end{align*}

By the Lipschitz condition of $\sigma$, \eqref{eq Z e estimate} and Cauchy-Schwarz's inequality, we have
\begin{align}\label{eq A e 21}
 \ee|A_{2,1}^\e(t,x)|^2 \le & \e \lambda^2(\e) L^2 \left( \int_0^t \sup_{(r,y)\in[0,s]\times\rr^d}\ee| Z^{\e,v^\e}(r,y) |^2\cdot \mathcal J(t-s)ds\right)\times\left(\int_0^t\|v^\e(s,\star)\|_{\HH}^2ds \right)^{\frac12}\notag\\
\le &  \e \lambda^2(\e)  C(L, N).
\end{align}

By the similar arguments as in the proof of \eqref{eq I 1}, we can show that
\begin{align}\label{eq A e 22}
\lim_{\e\rightarrow 0}\sup_{(t,x)\in[0,T]\times\rr^d}|A_{2,2}^{\e}(t,x)|\rightarrow0 \text{ in Probability.}
\end{align}

\vskip0.3cm
For the third term $A_3^{\e}$, by the Taylor's formula,  the Lipschitz continuity of $b'$ and \eqref{H3'}, we have
\begin{align*}
|A_3^\e(t,x)|\le &\Bigg|\int_0^t\Bigg\{\G_{\underline\alpha, \underline\delta}(t-s)\ast\Bigg[\frac{b\big(u^0(s,\star)+\sqrt\e\lambda(\e)Z^{\e,v^\e}(s,\star)\big) -b(u^0(s,\star) )}{\sqrt\e\lambda(\e)}\\&-b'(u^0(s,\star))Z^{\e,v^\e}(s,\star) \Bigg](x) \Bigg\}ds\Bigg|\notag\\
&+\Bigg|\int_0^t\left\{\G_{\underline\alpha, \underline\delta}(t-s)\ast\left[b'(u^0(s,\star))(Z^{\e,v^{\e}}(s,\star)-Z^v(s,\star))  \right](x) \right\}ds\Bigg|\notag\\
\le &L'\sqrt\e\lambda(\e)\int_0^t\G_{\underline\alpha, \underline\delta}(t-s)\ast  \left|Z^{\e,v^\e}(s,\star)\right|  (x) ds\notag\\
&+L\int_0^t \G_{\underline\alpha, \underline\delta}(t-s)\ast\left|Z^{\e,v^{\e}}(s,\star)-Z^v(s,\star)\right|(x) ds.
\end{align*}
Then,   by \eqref{eq Z e estimate}, we have
\begin{align}\label{eq A e 3}
\ee|A_3^\e(t,x)|^2
 \le & 2L'^2 \e\lambda^2(\e)  \int_0^t\int_{\rr^d}\G_{\underline\alpha, \underline\delta}(t-s,x-y) \sup_{(r,z)\in[0,s]\times\rr^d} \ee\left[\left|Z^{\e,v^\e}(r,z)\right|^2\right]   dsdy\notag\\
&+2L^2\int_0^t\int_{\rr^d} \G_{\underline\alpha, \underline\delta}(t-s,y)\sup_{(r,z)\in[0,s]\times\rr^d}\ee\left[\left|Z^{\e,v^{\e}}(r,z)-Z^v(r,z)\right|^2\right] dsdy\notag\\
\le & C(L') \e\lambda^2(\e) +2L^2\int_0^t\sup_{(r,z)\in[0,s]\times\rr^d}\ee\left[\left|Y^{\e,v^{\e},v}(r,z)\right|^2\right] ds.
\end{align}

Define the stopping time
$$
\tau^{M, \e}:=\inf\left\{t\le T; \sup_{(s,x)\in[0,t]\times\rr^d}|A_{2,2}^{\e}(s,x)|\ge M\right\},
$$
where $M$ is some constant large enough.

Putting \eqref{eq I0}, \eqref{eq A e 1}-\eqref{eq A e 3} together, we have
\begin{align*}
&\sup_{(s,x)\in[0,t]\times\rr^d}\ee\left[|Y^{\e,v^{\e},v}(s\wedge \tau^{M,\e},x)|^2\right]\ \notag\\
 \le&C\Bigg(\sup_{(s,x)\in[0,t]\times\rr^d}\ee\left[|A_1^\e(s\wedge \tau^{M,\e},x)|^2\right]+\sup_{(s,x)\in[0,t]\times\rr^d}\ee\left[|A_{2,1}^\e(s\wedge \tau^{M,\e},x)|^2\right]\notag\\
 &+\sup_{(s,x)\in[0,t ]\times\rr^d}\ee\left[|A_{2,2}^\e(s\wedge \tau^{M,\e},x)|^2\right]+\ee\left[|A_3^\e(t\wedge \tau^{M,\e},x)|^2\right]\Bigg)\notag\\
 \le & C\Bigg(  \lambda^{-2}(\e)+2\e \lambda^2(\e)+\sup_{(s,x)\in[0,t ]\times\rr^d}\ee\left[|A_{2,2}^\e(s\wedge \tau^{M,\e},x)|^2\right] \notag\\
 & +\int_0^t\sup_{(r,z)\in[0,s]\times\rr^d}\ee\left[\left|Y^{\e,v^{\e}v,}(r\wedge \tau^{M,\e},z)\right|^2ds\right] \Bigg).
\end{align*}
By Gronwall's lemma, we obtain that
\begin{align}\label{eq A e 4}
&\sup_{(s,x)\in[0,T]\times\rr^d}\ee\left[|Y^{\e,v^{\e},v}(s\wedge \tau^{M,\e},x) |^2\right]\notag\\
\le& C(T, L, L')\left(\lambda^{-2}(\e)+2\e \lambda^2(\e)+\sup_{(s,x)\in[0,T ]\times\rr^d}\ee\left[|A_{2,2}^\e(s\wedge \tau^{M,\e},x)|^2\right]\right).
\end{align}
Since $|A_{2,2}^\e(s\wedge \tau^{M,\e},x)|\le M $ for any $s\in[0,T]$, by \eqref{eq A e 22} and the  dominated convergence theorem, we know that
 \begin{align*}
 \sup_{(s,x)\in[0,T ]\times\rr^d}\ee\left[|A_{2,2}^\e(s\wedge \tau^{M,\e},x)|^2\right]
\longrightarrow 0,\ \ \text{as } \e\rightarrow 0.
 \end{align*}
This inequality, together with \eqref{eq A e 4}, implies that
\begin{align}
\sup_{(s,x)\in[0,T]\times\rr^d}\ee\left[|Y^{\e,v^{\e},v}(s\wedge \tau^{M,\e},x)|^2\right]\longrightarrow 0, \ \text{as } \e\rightarrow 0.
\end{align}
Letting $M\rightarrow\infty$, we obtain that for any $(t,x)\in[0,T]\times\rr^d$,
$$
\lim_{\e\rightarrow0}|Y^{\e,v^\e,v}(t,x)|= 0 \ \ \text{ in probability}.
$$

\vskip0.3cm

{\bf Step 2}. {\it Estimation of the increments}.

For any $ (t,x),(s,y)\in[0,T]\times K$ with $K$ being compact in $\rr^d$ and $t\ge s$,
\begin{align*}
&Y^{\e,v^\e,v}(t,x)-Y^{\e,v^\e,v}(s,y)\notag\\
=&\frac{1}{\lambda(\e)}\sum_{k\ge1}\Bigg\{\int_0^t\left\langle \G_{\underline\alpha, \underline\delta}(t-r,x-\star)\sigma\big(u^0(r,\star)+\sqrt\e\lambda(\e)Z^{\e,v^\e}(r,\star)\big), e_k  \right\rangle_{\HH}dB_r^k\notag\\
&\ \ \ \ \ \ \ \ \ -\int_0^s\left\langle \G_{\underline\alpha, \underline\delta}(s-r,y-\star)\sigma\big(u^0(r,\star)+\sqrt\e\lambda(\e)Z^{\e,v^\e}(r,\star)\big), e_k  \right\rangle_{\HH}dB_r^k\Bigg\}\notag\\
&+\Bigg\{\int_0^t\left\langle \G_{\underline\alpha, \underline\delta}(t-r,x-\star)\sigma\big(u^0(r,\star)+\sqrt\e\lambda(\e)Z^{\e,v^\e}(r,\star)\big), v^\e(r,\star)  \right\rangle_{\HH}dr\notag\\
&\ \ \ \  \ \ \ - \int_0^s\left\langle \G_{\underline\alpha, \underline\delta}(s-r,y-\star)\sigma\big(u^0(r,\star)+\sqrt\e\lambda(\e)Z^{\e,v^\e}(r,\star)\big), v^\e(r,\star)  \right\rangle_{\HH}dr\Bigg\}\notag\\
&-\Bigg\{\int_0^t\left\langle \G_{\underline\alpha, \underline\delta}(t-r,x-\star)\sigma\big(u^0(r,\star)\big), v(r,\star)  \right\rangle_{\HH}dr\notag\\
&\ \ \ \ \ \ \  -\int_0^s\left\langle \G_{\underline\alpha, \underline\delta}(s-r,y-\star)\sigma\big(u^0(r,\star)\big), v(r,\star)  \right\rangle_{\HH}ds\notag\Bigg\}\\
&+\Bigg\{\int_0^t\G_{\underline\alpha, \underline\delta}(t-r)\ast\left[\frac{b\big(u^0(r,\star)+\sqrt\e\lambda(\e)Z^{\e,v^\e}(r,\star)\big) -b(u^0(r,\star) )}{\sqrt\e\lambda(\e)} \right](x) dr\notag\\
&\ \ \ \ \ \ \ \ -\int_0^s\G_{\underline\alpha, \underline\delta}(s-r)\ast\left[\frac{b\big(u^0(r,\star)+\sqrt\e\lambda(\e)Z^{\e,v^\e}(r,\star)\big) -b(u^0(r,\star) )}{\sqrt\e\lambda(\e)} \right](y) dr\Bigg\}\notag\\
&-\Bigg\{\int_0^t\G_{\underline\alpha, \underline\delta}(t-r)\ast\left[b'(u^0(r,\star))Z^v(r,\star)  \right](x) dr \notag\\
&\ \ \ \ \ \ -\int_0^s\G_{\underline\alpha, \underline\delta}(s-r)\ast\left[b'(u^0(r,\star))Z^v(r,\star)  \right](y) dr\Bigg\}\notag\\
=:& \frac1{\lambda(\e)}B_1^\e(t,s,x,y)+B_2^\e(t,s,x,y)+B_3^\e(t,s,x,y)+B_4^\e(t,s,x,y)+B_5^\e(t,s,x,y).
\end{align*}

For the first term, by the Lipschitz continuity of $\sigma$, \eqref{eq Z e estimate} and Lemma \ref{lem holder},   we obtain that for any $\beta_1\in ]0,(1-\eta)/2[, 0<\beta_2<\min\{(1-\eta)\alpha_0/2, 1/2\}$
\begin{align}\label{eq B 1}
&\ee\left[|B_1^\e(t,s,x,y)|^2\right]\notag\\
\le &C_1\ee\Bigg[\int_0^s \|(\G_{\underline\alpha, \underline\delta}(t-r,x-\star)-\G_{\underline\alpha, \underline\delta}(s-r,x-\star))\sigma\big(u^0(r,\star)+\sqrt\e\lambda(\e)Z^{\e,v^\e}(r,\star)\big)\|_{\HH}^2dr\notag\\
 &+\int_s^t \|\G_{\underline\alpha, \underline\delta}(t-r,x-\star) \sigma\big(u^0(r,\star)+\sqrt\e\lambda(\e)Z^{\e,v^\e}(r,\star)\big)\|_{\HH}^2dr\notag\\
&+\int_0^s \|(\G_{\underline\alpha, \underline\delta}(s-r,x-\star)-\G_{\underline\alpha, \underline\delta}(s-r,y-\star))\sigma\big(u^0(r,\star)+\sqrt\e\lambda(\e)Z^{\e,v^\e}(r,\star)\big)\|_{\HH}^2dr\Bigg]\notag\\
\le & C_2\sup_{(r,z)\in[0,T]\times\rr^d}\ee\left[|\sigma (u^0(r,\star)+\sqrt\e\lambda(\e)Z^{\e,v^\e}(r,\star) )|^2\right]\notag\\
&\times\Bigg[\int_0^s \|(\G_{\underline\alpha, \underline\delta}(t-r,x-\star)-\G_{\underline\alpha, \underline\delta}(s-r,x-\star))\|_{\HH}^2dr+\int_s^t \|\G_{\underline\alpha, \underline\delta}(t-r,x-\star) \|_{\HH}^2dr\notag\\
&\ \ \ \ +\int_0^s \|(\G_{\underline\alpha, \underline\delta}(s-r,x-\star)-\G_{\underline\alpha, \underline\delta}(s-r,y-\star))\|_{\HH}^2dr\Bigg]  \notag\\
\le & C_3(|t-s|^{2\beta_1}+|x-y|^{2\beta_2}),
\end{align}
where $C_3$ is independent of $\e>0$.

Using the similar argument in \eqref{eq B 1} and noticing the fact $v^\e,v\in\AA_{\HH}^N$, we give the following estimates without the proof
\begin{align}\label{eq B 2}
\ee\left[|B_i^\e(t,s,x,y)|^2\right]
\le  C(|t-s|^{2\beta_1}+|x-y|^{2\beta_2}),\ \ i=2,3,
\end{align}
where $C$ is independent of $\e>0$.

We also can deal with the last two terms $B_4^{\e}, B_5^\e$ by the same argument. Next we only estimate $B_4^\e$, and  the same result holds for $B_5^\e$.

By Taylor's formula,  \eqref{H3'} and
using the same approach in estimation of $\Lambda_3$ in  the proof of \cite[Proposition  2.10]{EM}, we have for any $t\in[0,T]$,
\begin{align}\label{eq B 4}
&\ee\left[|B_4^\e(t,s,x,y)|^2\right]\notag\\
\le&  L\ee\Bigg|\int_0^t\G_{\underline\alpha, \underline\delta}(t-r)\ast|Z^{\e,v^\e}(r,\star)| (x) dr -\int_0^s\G_{\underline\alpha, \underline\delta}(s-r)\ast|Z^{\e,v^\e}(r,\star)|(y) dr\Bigg|^2\notag\\
\le& C\left[|t-s|+\int_0^t\sup_{z\in\rr^d}\ee|Z^{\e,v^\e}(t-s+r,x-y+z)-Z^{\e,v^\e}(r,z)|^2dr\right]\notag\\
\le & C\left[ |t-s|+T(|t-s|^{2\beta_1}+|x-y|^{2\beta_2})\right],
\end{align}
where in the last inequality we have used the estimate in the proof of Theorem 3.1 in \cite{BEM}, $C$ is independent of $\e$.

  Putting together all the estimates, we get \eqref{eq Z 2}.
The proof is complete.
\nprf

\section{ Appendix }
   To make reading easier, we present here some results on the kernel $G$    from  Boulanba {\it et al.} \cite{BEM}.

\begin{lemma}\cite[Lemma 1.1]{BEM}\label{lem fundamental solution} For $\alpha\in ]0,2]\setminus\{1\}$ such that $|\delta|\le \min\{\alpha,2-\alpha\}$, the following statements hold.
\begin{itemize}
  \item[(i)] The function $G_{\alpha,\delta}(t,x)$ is the density of a L\'evy $\alpha$-stable process in time $t$.
  \item[(ii)]
Semigroup property: $G_{\alpha,\delta}(t,x)$  satisfies the Chapman-Kolmogorov equation, i.e., for $0<s<t$,
  $$
  G_{\alpha,\delta}(t+s,x)=\int_{\rr}G_{\alpha,\delta}(t,y)G_{\alpha,\delta}(s,y-x)dy.
  $$
  \item[(iii)] Scaling property: $G_{\alpha,\delta}(t,x)=t^{-1/\alpha}G_{\alpha,\delta}(1,t^{-1/\alpha}x)$.
  \item[(iv)] There exists a constant $c_{\alpha}$ such that $0\le G_{\alpha,\delta}(1,x)\le c_{\alpha}/(1+|x|^{1+\alpha})$, for all $x\in\rr$.
\end{itemize}

\end{lemma}

The next proposition studies the H\"older regularity of the Green function, whose proof is contained in  the proof of Proposition 3.2 in \cite{BEM}.
\begin{lemma}\label{lem holder} Under $(H_{\eta}^{\underline{\alpha}})$, it holds that
\begin{itemize}
  \item[(i)] For each $0\le s<t\le T,\beta_1\in ]0,(1-\eta)/2[$, there exists a constants $C>0$ such that
  $$
  \int_0^s \| \G_{\underline\alpha, \underline \delta}(t-r, \star)-\G_{\underline\alpha, \underline \delta}(s-r, \star)\|_{\HH}^2dr\le C |t-s|^{2\beta_1},
  $$
  and
  $$
  \int_s^{t} \| \G_{\underline\alpha, \underline \delta}(t-r, \star)\|_{\HH}^2dr\le C |t-s|^{2\beta_1}.
  $$
  \item[(ii)] For each $0<\beta_2<\min\{(1-\eta)\alpha_0/2, 1/2\}$, there exists a constant $C>0$ for any $x,y\in\rr^d$,
  $$
  \int_0^T \| \G_{\underline\alpha, \underline \delta}(T-s, x-\star)-\G_{\underline\alpha, \underline \delta}(T-s, y-\star)\|_{\HH}^2ds\le C |x-y|^{2\beta_2}.
  $$

\end{itemize}

\end{lemma}

The next lemma is about the H\"older regularity of the stochastic integral. See \cite[Proposition 3.2]{BEM}. For a given predictable random field $V$, we set
$$
U(t,x):=\sum_{k\ge1}\int_0^t \langle\G_{\underline{\alpha},\underline{\delta}}(t-s,x-\star)V(s,\star),e_k\rangle_{\HH}dB_s^k.
$$
\begin{lemma}\label{lem Holder stochastic} Assume that $\sup_{0\le t\le T}\sup_{x\in\rr^d}\ee(|V(t,x)|^p)$ is finite for some $p$ large enough. Then under $(H_{\eta}^{\underline{\alpha}})$, we have
\begin{itemize}
  \item[(i)] For each $x\in\rr^d$ a.s. the mapping $t\rightarrow U(t,x)$ is  $\beta_1$-H\"older continuous for $0<\beta_1<(1-\eta)/2$.
  \item[(ii)] For each $t\in[0,T]$ a.s. the mapping $x\rightarrow U(t,x)$ is $\beta_2$-H\"older continuous for $0<\beta_2<\min\{\alpha_0(1-\eta)/2,1/2\}$.
\end{itemize}

\end{lemma}

The following  result  is a consequence of   Lemma A.1 in  \cite{BMS}.
\begin{lemma}\label{Lem 3}\rm{
Let $K$ be a compact set in $\rr^d$ and let $\{V^\e(t,x):(t,x)\in[0,T]\times K\}$ be a family of real-valued functions. Assume
\begin{itemize}
  \item[(A1).] For any $(t,x)\in[0,T]\times K$,
     $$\lim_{\e\rightarrow 0}|V^{\e}(t,x)|=0.$$
  \item[(A2).] There exist $\beta_1,\beta_2>0$ satisfied that for any $(t,x),(t',x')\in[0,T]\times K$,
$$
|V^\e(t,x)-V^\e(t',x')|\le C(|t-t'|^{\beta_1}+|x-x'|^{\beta_2}),
$$
where $C$ is a constant independent of $\e$.
\end{itemize}
Then for any $\theta\in(0,1)$, we have
$$\lim_{\e\to0}\|V^\e\|_{\theta\beta_1,\theta\beta_2}=0. $$
}\end{lemma}
\bprf
The stochastic version of this   lemma in  \cite[Lemma A1]{BMS} is  proved by  the Garsia-Rodemich-Rumsey's lemma.  Here we give  a direct proof for the deterministic case.

In view of Arzel\`a-Ascoli theorem, a sequence in $C([0,T]\times K;\rr) $ converges uniformly if and only if it is equicontinuous and converges pointwise.  Thence, under conditions (A.1) and (A.2), we know that
$$\lim_{\e\rightarrow0}\sup_{(t,x)\in[0,T]\times K}|V^\e(t,x)|\rightarrow0.
$$
Thus, for any $\theta\in(0,1),(t,x)\neq(t',x')\in[0,T]\times K$, we have
\begin{align*}
\frac{ |V^\e(t,x)-V^\e(t',x')|}{(|t-t'|^{\beta_1}+|x-x'|^{\beta_2}) ^{\theta}}= &\frac{|V^\e(t,x)-V^\e(t',x')|^{\theta} |V^\e(t,x)-V^\e(t',x')|^{1-\theta}}{(|t-t'|^{\beta_1}+|x-x'|^{\beta_2}) ^{\theta}}\\
 \le &C  (|V^\e(t,x)|+|V^\e(t',x')|)^{1-\theta}\\
 &\longrightarrow0 \ \ \ \ \ \ \text{uniformly on } ([0,T]\times K)^2.
\end{align*}
The proof is complete.
\nprf

\vskip0.3cm
\noindent{\bf Acknowledgements:}
  Y. Li and S. Zhang were supported by Natural Science Foundation of China (11471304, 11401556). R. Wang was supported by Natural Science Foundation of China (11301498, 11431014) and the Fundamental Research Funds for the Central Universities (WK 0010000048). N. Yao was supported by Natural Science Foundation of China (11101313).

\end{document}